\documentclass{article}
\usepackage{amsmath,amssymb,amsthm,amsfonts,latexsym}
\usepackage[all]{xy}
\usepackage{amscd,graphicx,stmaryrd}

\begin{document}

\newcommand{\nc}[2]{\newcommand{#1}{#2}}
\newcommand{\rnc}[2]{\renewcommand{#1}{#2}}
\def\note#1{{}}
\def\Label#1{\label{#1}\ifmmode\llap{[#1] }\else
\marginpar{\smash{\hbox{[#1]}}}\fi}

\newtheorem{de}{Definition}[section]
\newtheorem{prop}[de]{Proposition}
\newtheorem{lem}[de]{Lemma}
\newtheorem{corollary}[de]{Corollary}
\newtheorem{them}[de]{Theorem}
\newtheorem{example}[de]{\it Example}
\newtheorem{rem}[de]{Remark}


\rnc{\[}{\begin{equation}}
\rnc{\]}{\end{equation}}
\nc{\wegengruen}{\end{equation}}

\nc{\be}{\begin{enumerate}}
\nc{\ee}{\end{enumerate}}
\nc{\bd}{\begin{diagram}}
\nc{\ed}{\end{diagram}}
\nc{\bi}{\begin{itemize}}
\nc{\ei}{\end{itemize}}
\nc{\bpr}{\begin{prop}}
\nc{\bth}{\begin{them}}
\nc{\ble}{\begin{lem}}
\nc{\bco}{\begin{corollary}}
\nc{\bre}{\begin{remark}}
\nc{\bex}{\begin{example}}
\nc{\bde}{\begin{de}}
\nc{\ede}{\end{de}}
\nc{\epr}{\end{prop}}
\nc{\ethe}{\end{them}}
\nc{\ele}{\end{lem}}
\nc{\eco}{\end{corollary}}
\nc{\ere}{\hfill\mbox{$\losenge$}\end{remark}}
\nc{\eex}{\hfill\mbox{$\losenge$}\end{example}}
\nc{\bpf}{{\it Proof.~~}}
\nc{\epf}{\hfill\mbox{$\square$}\vspace*{3mm}}
\nc{\hsp}{\hspace*}
\nc{\vsp}{\vspace*}

\nc{\ot}{\otimes}
\nc{\te}{\!\ot\!}
\nc{\bmlp}{\mbox{\boldmath$\left(\right.$}}
\nc{\bmrp}{\mbox{\boldmath$\left.\right)$}}
\nc{\LAblp}{\mbox{\LARGE\boldmath$($}}
\nc{\LAbrp}{\mbox{\LARGE\boldmath$)$}}
\nc{\Lblp}{\mbox{\Large\boldmath$($}}
\nc{\Lbrp}{\mbox{\Large\boldmath$)$}}
\nc{\lblp}{\mbox{\large\boldmath$($}}
\nc{\lbrp}{\mbox{\large\boldmath$)$}}
\nc{\blp}{\mbox{\boldmath$($}}
\nc{\brp}{\mbox{\boldmath$)$}}
\nc{\LAlp}{\mbox{\LARGE $($}}
\nc{\LArp}{\mbox{\LARGE $)$}}
\nc{\Llp}{\mbox{\Large $($}}
\nc{\Lrp}{\mbox{\Large $)$}}
\nc{\llp}{\mbox{\large $($}}
\nc{\lrp}{\mbox{\large $)$}}
\nc{\lbc}{\mbox{\Large\boldmath$,$}}
\nc{\lc}{\mbox{\Large$,$}}
\nc{\Lall}{\mbox{\Large$\forall\;$}}
\nc{\bc}{\mbox{\boldmath$,$}}
\nc{\ra}{\rightarrow}
\nc{\ci}{\circ}
\nc{\cc}{\!\ci\!}
\nc{\lra}{\longrightarrow}
\nc{\imp}{\Rightarrow}
\rnc{\iff}{\Leftrightarrow}
\nc{\inc}{\mbox{$\,\subseteq\;$}}
\rnc{\subset}{\inc}
\nc{\0}{\sp{(0)}}
\nc{\3}{\sp{(3)}}
\nc{\4}{\sp{(4)}}
\nc{\5}{\sp{(5)}}
\nc{\6}{\sp{(6)}}
\nc{\7}{\sp{(7)}}
\newcommand{\Boxneu}{\square}
\def\tr{{\rm tr}}
\def\Tr{{\rm Tr}}
\def\st{\stackrel}
\def\<{\langle}
\def\>{\rangle}
\def\d{\mbox{$\mathop{\mbox{\rm d}}$}}
\def\pr{\mbox{$\mathop{\mbox{\rm pr}}$}}
\def\id{\mbox{$\mathop{\mbox{\rm id}}$}}
\def\ker{\mbox{$\mathop{\mbox{\rm Ker$\,$}}$}}
\def\hom{\mbox{$\mathop{\mbox{\rm Hom}}$}}
\def\corep{\mbox{$\mathop{\mbox{\rm Corep}}$}}
\def\Mat{\mathrm{Mat}}
\def\im{\mathrm{i}}
\def\map{\mbox{$\mathop{\mbox{\rm Map}}$}}
\def\spec{\mbox{$\mathop{\mbox{\rm spec}}$}}
\def\o{\sp{[1]}}
\def\t{\sp{[2]}}
\def\mo{\sp{[-1]}}
\def\z{\sp{[0]}}

\nc{\spp}{\mbox{${\cal S}{\cal P}(P)$}}
\nc{\ob}{\mbox{$\Omega\sp{1}\! (\! B)$}}
\nc{\op}{\mbox{$\Omega\sp{1}\! (\! P)$}}
\nc{\oa}{\mbox{$\Omega\sp{1}\! (\! A)$}}
\nc{\dr}{\mbox{$\Delta_{R}$}}
\nc{\dsr}{\mbox{$\Delta_{\Omega^1P}$}}
\nc{\ad}{\mbox{$\mathop{\mbox{\rm Ad}}_R$}}
\nc{\as}{\mbox{$A(S^3\sb s)$}}
\nc{\bs}{\mbox{$A(S^2\sb s)$}}
\nc{\slc}{\mbox{$A(SL(2,\C))$}}
\nc{\SUq}{\mbox{$\cO(\mathrm{SU}_q(2))$}}
\nc{\CSU}{\mbox{$\cC(\mathrm{SU}_q(2))$}}
\nc{\podl}{\mbox{$\cO(\mathrm{S}^2_{qs})$}}
\nc{\Dq}{\mbox{$\cO(\mathrm{D}^2_{q})$}}
\nc{\Cpodl}{\mbox{$\cC(\mathrm{S}^2_{qs})$}}
\nc{\CS}{\mbox{$\cC(\mathrm{S}^1)$}}
\nc{\tc}{\widetilde{can}}
\def\slq{\mbox{$\cO(SL_q(2))$}}
\def\asq{\mbox{$\cO(S_{q,s}^2)$}}
\def\esl{{\mbox{$E\sb{\frak s\frak l (2,{\Bbb C})}$}}}
\def\esu{{\mbox{$E\sb{\frak s\frak u(2)}$}}}
\def\ox{{\mbox{$\Omega\sp 1\sb{\frak M}X$}}}
\def\oxh{{\mbox{$\Omega\sp 1\sb{\frak M-hor}X$}}}
\def\oxs{{\mbox{$\Omega\sp 1\sb{\frak M-shor}X$}}}
\def\Fr{\mbox{Fr}}

\rnc{\epsilon}{\varepsilon}
\rnc{\phi}{\varphi}
\nc{\ha}{\mbox{$\alpha$}}
\nc{\hb}{\mbox{$\beta$}}
\nc{\hg}{\mbox{$\gamma$}}
\nc{\hc}{\mbox{$\gamma$}}
\nc{\hd}{\mbox{$\delta$}}
\nc{\he}{\mbox{$\varepsilon$}}
\nc{\hz}{\mbox{$\zeta$}}
\nc{\hs}{\mbox{$\sigma$}}
\nc{\hk}{\mbox{$\kappa$}}
\nc{\hm}{\mbox{$\mu$}}
\nc{\hn}{\mbox{$\nu$}}
\nc{\hl}{\mbox{$\lambda$}}
\nc{\hG}{\mbox{$\Gamma$}}
\nc{\hD}{\mbox{$\Delta$}}
\nc{\hT}{\mbox{$\Theta$}}
\nc{\ho}{\mbox{$\omega$}}
\nc{\hO}{\mbox{$\Omega$}}
\nc{\hp}{\mbox{$\pi$}}
\nc{\hP}{\mbox{$\Pi$}}

\def\C{{\mathbb C}}
\def\N{{\mathbb N}}
\def\R{{\mathbb R}}
\def\Z{{\mathbb Z}}
\def\T{{\mathbb T}}
\def\cO{{\mathcal O}}
\def\A{{\mathcal A}}
\def\O{\cO}
\def\cC{{\mathcal C}}
\def\T{{\cal T}}
\def\S{{\mathfrak S}}
\def\cK{{\cal K}}
\def\cH{{\cal H}}\def\ch{{\cal H}}
\def\H{{\cal H}}
\def\ta{\tilde a}
\def\tb{\tilde b}
\def\td{\tilde d}

\nc{\Cs}{\mbox{$\mathrm{C}^*$}}
\nc{\lN}{\mbox{${\ell}^2(\N_0)$}}
\nc{\lZ}{\mbox{${\ell}^2(\Z)$}}
\nc{\mpt}{\hspace{-1pt}}
\nc{\ppt}{\hspace{1pt}}
\nc{\lin}{\mbox{$\mathrm{span}$}}
\nc{\trans}{\mathrm{t}}
\nc{\pink}{\mbox{$\Pi_{k=0}^{n-1}$}}
\nc{\Kn}{\mbox{$\mathrm{K}_0$}}
\nc{\KN}{\mbox{$\mathrm{K}^0$}}
\nc{\ovM}{\overline{M}}


\title{Fibre product approach to index pairings for the generic Hopf fibration of
$\mathrm{SU}_q(2)$}

\author{Elmar Wagner}

\date{}

\maketitle

\thispagestyle{empty}

\mbox{ }\\[-24pt]
\centerline{{\small Instituto de F\'isica y Matem\'aticas}}
\centerline{{\small Universidad Michoacana de San Nicol\'as de Hidalgo, Morelia, M\'exico }}
\centerline{{\small E-mail: Elmar.Wagner@math.uni-leipzig.de}\\}


\begin{abstract}
A fibre product construction is used to 
give a description of quantum line bundles over 
the generic Podle\'s spheres by gluing two quantum discs along their boundaries. 
Representatives of the corresponding $K_0$-classes are given in terms of 
1-dimensional projections belonging to the C*-algebra, and in terms of analogues 
of the classical Bott projections. The $K_0$-classes of 
quantum line bundles derived from the generic Hopf fibration of quantum 
$\mathrm{SU}(2)$ are determined and the index pairing is computed. 
It is argued that taking the projections 
obtained from the fibre product construction 
yields a significant simplification of earlier index computations. 
\end{abstract}

\noindent
{\small 
}

\section{Introduction}

The main goal of the paper is to show that the fibre product approach 
provides an effective tool for simplifying index computations. 
Our example is one of the most extensively studied quantum spaces, 
namely the Podle\'s spheres \cite{Podles}. 
These two parameter deformations of the classical 2-sphere served as a 
guiding example in the development of the theory of 
coalgebra bundles (quantum principal bundles) \cite{Brz,BM}. 
In this setting, the Podle\'s spheres  and the quantum $\mathrm{SU}(2)$ 
play the roles of base and total space of a noncommutative Hopf fibration, 
respectively. A projective module description of all 
associated quantum line bundles was given in \cite{HM} for 
the standard Podle\'s sphere, and in \cite{HMS,SW} for general parameters. 
Since the cyclic cohomology of the Podle\'s spheres was known \cite{MNW}, 
the projective module description made it possible to compute 
the index pairing (Chern-Connes pairing) between cyclic cohomology and 
the $\Kn$-classes of the quantum line bundles. 
This was done first in \cite{H} for the standard Podle\'s sphere, and then 
extended to the general case in \cite{HMS}. Their proofs 
relied heavily on the index theorem, arguing that 
the index is an integer and taking 
a suitable limit of the parameters. 
The main obstacle for a direct computation was the growing size 
of the projection matrices, leading to complicated expressions 
which are difficult to handle. 
It is therefore desirable to find a tool for obtaining simpler 
representatives of $K$-theory classes. 

The computations of the $K$-groups of quantum 2-spheres in 
\cite{HMS06} and of quantum 3-spheres in \cite{BHMS}  
strongly encourage to use a fibre product approach 
to index pairings. 
Furthermore, an 
explicit description of projective modules 
and a significant simplification of
the index problem were obtained in \cite{DHHM}
by using Bass' construction of 
the Mayer-Vietoris boundary map in (algebraic) 
$K$-theory \cite{Ba}. 
In the present paper, we shall see another 
example where the fibre product approach simplifies  
the index computation considerably. 

The fact that the C*-algebra of generic Podle\'s spheres is 
given by a fibre product of two Toeplitz algebras 
follows from \cite[Proposition 1.2]{s-a91}. 
In Section \ref{sec-glue}, we explain how this description 
can be interpreted as gluing two quantum discs 
along the boundary (see also \cite{bk96,cm00}). 
Our first result concerns the construction 
of quantum line bundles. Since the boundaries of the 
quantum discs are identified with the 
classical circle $\mathrm{S}^1$, 
we can glue two trivial bundles  
over the quantum disc 
by using the same transition function 
as in the commutative case. 
The gluing procedure is described by a fibre product 
and the resulting projective modules are considered as 
quantum line bundles of winding number $N$, where 
$N\in\Z$ is the degree of the transition function. 
We then show that the quantum line bundles 
are isomorphic to projective modules 
given by elementary 1-dimensional projections. 
In a sense, the quantum case is simpler 
than its classical counterpart since there are no 
non-trivial 1-dimensional projections 
in $\cC(\mathrm{S}^2)$. 

The link between the fibre product approach 
of quantum line bundles and the Hopf fibration 
of $\mathrm{SU}_q(2)$ will be established in Section \ref{pmd}. 
We prove that, for each winding number, the 
projection derived from the noncommutative Hopf fibration 
is Murray-von Neumann equivalent to the 
1-dimensional projection of the fibre product construction. 
As an intermediate step, we introduce Murray-von Neumann equivalent 
$2\times 2$-projections whose entries are rational functions 
of the generators of the Podle\'s spheres. These 
$2\times 2$-projections can be 
viewed as analogues of the classical Bott projections. 
Moreover, their explicit description enables us to give a 
short direct proof of the statement that 
the direct sum of quantum line bundles with winding number 
$N$ and $-N$ is equivalent to a free rank 2 module. 

In the final section, we carry out the index computations. 
Because of the previous results, we are free to choose the 
most convenient representatives of the $\Kn$-classes. By using the 
1-dimensional projections from the fibre product approach, 
the index pairing reduces to its simplest possible form;  
it remains to calculate a trace of a projection onto a finite 
dimensional subspace. In analogy to the classical case, 
one generator of $K$-homology computes the winding 
number, and the other detects the rank of the 
projective module.

\section{Preliminaries}
                                                             \label{prel}

\subsection{Fibre products}                    \label{sec-fp}
Let $\pi_0: A_0 \rightarrow A_{01}$ and  $\pi_1:A_1\rightarrow
A_{01}$ be morphisms of C*-algebras. The fibre product
$A:=A_0{\times}_{(\pi_0,\pi_1)} A_1$ is defined by the
pull-back diagram

\begin{equation}                                   \label{A_is_fibre_product}
        \begin{CD}
    {A} @ >{\mathrm{pr}_1}>> {A_1} @.\\
    @ V{\mathrm{pr}_0} VV @ V{\pi_1} VV @.\\
    {A_0} @ >{\pi_0} >> {A_{01}}. @.\\
        \end{CD}
\end{equation}
Up to a unique isomorphism, $A$ is given by
\begin{equation}                                       \label{A}
A=\left\lbrace (a_0,a_1)\in A_0\times A_1 : \pi_0(a_0)=\pi_1(a_1)
\right\rbrace,
\end{equation}
where $A_0\times A_1$ denotes the outer direct sum of C*-al\-ge\-bras
with multiplication $(a_0,a_1)(b_0,b_1)=(a_0b_0,a_1b_1)$
and involution $(a_0,a_1)^{*}=(a_0^{*},a_1^{*})$.
The morphisms
$\mathrm{pr}_0:A\rightarrow A_0$ and
$\mathrm{pr}_1:A\rightarrow A_1$ are then the left and right projections,
respectively.

Similarly, if $B$ is a C*-algebra, and
$\pi_0: A_0 \rightarrow A_{01}$ and  $\pi_1:A_1\rightarrow A_{01}$
are morphisms of left $B$-modules, then the fibre product
of vector spaces
$A:=A_0{\times}_{(\pi_0,\pi_1)} A_1$
is a left $B$-module
with left action given by $b.(a_0,a_1)=(b. a_0,b. a_1)$, where
the dot denotes the left action.

Now let $A_j$ be left $B_j$-modules,
$j=0,1,01$, and let $\rho_0: B_0 \rightarrow B_{01}$ and  $\rho_1:B_1\rightarrow
B_{01}$ be morphisms of C*-algebras.
Suppose that the vector space morphisms
$\pi_0: A_0 \rightarrow A_{01}$ and  $\pi_1:A_1\rightarrow A_{01}$
intertwine the actions, that is,
$\pi_i(b_i.a_i)=\rho_i(b_i).\pi_i(a_i)$ for all $b_i\in B_i$ and $a_i\in A_i$, $i=0,1$.
Then $A$ is a left
$B=B_0{\times}_{(\rho_0,\rho_1)} B_1$-module. This follows from the preceding
with the left $B$-actions given by
$b.a_i:= \mathrm{pr}_i(b).a_i$ for $a_i\in\A_i$, $i=0,1$, and
$b.a_{01}:=\rho_0\!\circ\!\mathrm{pr}_0(b).a_{01}=\rho_1\!\circ\!\mathrm{pr}_1(b).a_{01}$
for $a_{01}\in\A_{01}$.

\subsection{The generic Podle\'s spheres}                           \label{Ps}

Let $q\in(0,1)$ and  $s\in [0,\infty)$.
The *-algebra $\podl$ of polynomial functions on the
quantum 2-sphere $\mathrm{S}^2_{qs}$ is
generated by $\eta_s$, $\eta_s^*$
and an hermitian element $\zeta_s$ satisfying the relations \cite{Podles}
\[                                                      \label{rel}
\zeta_s \eta_s = q^2 \eta_s \zeta_s,\quad \!\!
\eta_s^*\eta_s=(1-\zeta_s)(s^2+\zeta_s),\quad \!\!
\eta_s\eta_s^*=(1-q^{-2}\zeta_s)(s^2+q^{-2}\zeta_s).
\]
The case $s=0$ corresponds to the standard Podle\'s sphere.
Its algebraic and \Cs-algebraic properties differ considerably
from the other Podle\'s spheres which we call ``generic''.
The standard Podle\'s sphere is treated in \cite{HW} and will be excluded here.
Note that, for $s>0$, $\podl$ and
$\cO(\mathrm{S}^2_{q,s^{-1}})$ are isomorphic,
where an isomorphism is
given by $\zeta_s\rightarrow -s^{2} \zeta_{s^{-1}}$ and
$\eta_s\rightarrow s\ppt\eta_{s^{-1}}$. For practical reasons,
we shall consider both  versions of the same *-algebra.

Recall that the universal \Cs-algebra of some *-algebra is defined as
the \Cs-completion with respect to the
universal \Cs-norm given by the
supremum of the operator norms over all bounded
*-re\-pre\-sen\-ta\-tions (if the supremum exists).
The universal \Cs-algebra of $\podl$ will be denoted by $\Cpodl$.

The *-re\-pre\-sen\-ta\-tions of $\podl$
were described in \cite{MNW,Podles} and are classified as follows.
Given a Hilbert space  $\H_0$,  set
$\H:=\oplus_{n=0}^\infty \H_n$, where $ \H_n=\H_0$.
For $e\in \H_0\setminus\{0\}$ and $n\in\N_0$, we denote by $e_n$ the vector in
$\H$ which has $e$ as its $n$-th component and zero otherwise.
On $\H$, the shift operator $S$ is
given by
\begin{equation}                             \label{S}
 Se_n=e_{n+1}.
\end{equation}

Now let $\H_0^{+}$, $\H_0^{-}$, $\H^0$ be Hilbert spaces and
let $U$ be a unitary operator on $\H^0$. Set
$\H^{+}:= \oplus_{n=0}^\infty \H_n^{+}$ and
$\H^{-}:= \oplus_{n=0}^\infty \H_n^{-}$.
For $s\neq 0$, a *-re\-pre\-sen\-ta\-tion $\rho$ of $\podl$
on a Hilbert space $\H$ has
the direct sum decomposition  $\rho=\rho_-\oplus\rho_0\oplus\rho_+$
on $\H=\H^-\oplus\H^0\oplus \H^+$ and is given by
\begin{align}                           \nonumber              
&\rho_-(\zeta_s)e^-_n =  -s^2 q^{2(n+1)} e^-_n, \quad
\rho_-(\eta_s)e^-_n= s \sqrt{(1+s^2q^{2(n+1)} )(1-q^{2(n+1)})}\,S e^-_{n},\\
&\rho_0(\zeta_s)e^0= 0,\quad \rho_0(\eta_s)e^0= sU, \nonumber\\
&\rho_+(\zeta_s)e^+_n=  q^{2(n+1)} e^+_n, \quad      \nonumber  
\rho_+(\eta_s)e^+_n= \sqrt{(1-q^{2(n+1)} )(s^2+q^{2(n+1)})}\,  S e^+_{n}.
\end{align}
The representations $\rho_-$, $\rho_0$ and $\rho_+$
are irreducible if and only if the Hilbert spaces
$\H^+_0$, $\H^-_0$ and $\H^0$ are 1-dimensional.

It was argued in \cite{MNW,s-a91} that the direct sum
$\rho_-\oplus\rho_+$ of the irreducible representations
$\rho_-$ and $\rho_+$ on $\lN$ yields a faithful
representation of $\Cpodl$. Therefore
we can consider $\Cpodl$ as the concrete \Cs-algebra
given by the representation $\rho_-\oplus\rho_+$
on $\lN\oplus \lN$. By a slight abuse of notation, we shall
frequently identify elements of $\Cpodl$ with its image
under the mapping $\rho_-\oplus\rho_+$.

For later use, let us also mention that
\[                                                          \label{spec}
\spec(\zeta_s)= \{-s^2 q^2, -s^2 q^4, \ldots \}\cup \{0\}\cup
\{q^2, q^4, \ldots \}.
\]

To give an explicit description of $\Cpodl$, we make use of
the quantum disc.
The *-algebra $\Dq$ of polynomial functions on the quantum disc
is generated by two generators $z$ and $z^{*}$ with relation
\[                                                          \label{Dq}
z^* z -q^2 z z^* = 1- q^2.
\]
The universal \Cs-algebra of $\Dq$
is well known. It has been discussed by several authors
(see, e.g., \cite{kl93,MNW,s-a91}) that it
is isomorphic to the Toeplitz
algebra $\T$, where $\T$ can be characterized as
the universal \Cs-algebra generated by the
unilateral shift $S$ on $\lN$.

Let $U$ denote the unitary generator of $\CS$.
The so-called symbol map is the  *-homomorphism
$\hs:\T\rightarrow\CS$
defined by $\hs(S)=U$.
In \cite{s-a91},
it has been shown that (for $s>0$)
\[
\Cpodl\cong                                                  \label{CS}
\{(a_0,a_1)\in \T\oplus \T \,:\, \hs(a_0)=\hs(a_1)\}.
\]
The *-re\-pre\-sen\-ta\-tions $\rho_-$ and $\rho_+$ of $\Cpodl$ on $\lN$ are
then given by
\begin{equation}                                  \label{rhopm}
 \rho_-((a_0,a_1))=a_0,\quad \rho_+((a_0,a_1))=a_1, \quad
 (a_0,a_1)\in\Cpodl.
\end{equation}

The K-theory and K-homology of $\Cpodl$ has been computed in \cite{MNW}.
There it is shown that $\Kn(\Cpodl)\cong \Z \oplus \Z$
and $\KN(\Cpodl)\cong \Z \oplus \Z$.  For  $s>0$,
the two  generators of the $\Kn$-group are given by the class $[1]$ of
the unit element $1\in\Cpodl$, and by the class  $[(0,1-SS^{*})]$ of
the 1-dimensional projection onto $\C(0,e_0)\subset\lN\oplus\lN$.

Describing an even Fredholm module by a pair of representations
on the same Hilbert space,
one generator of $\KN(\Cpodl)$ is given
by $[(\rho_+,\rho_-)]$ on $\lN$. The second generator is
obtained by a pull-back of an even Fredholm module
on $\CS$ via the symbol map
$\hs:\Cpodl\rightarrow \CS$,\ \,$\hs((a_0,a_1)):=\hs(a_0)=\hs(a_1)$.
The even  Fredholm module on $\CS$ was described in
\cite{MNW1} by the following pair of representations
$\pi_\pm:\CS\rightarrow \lZ$:
\begin{align*}
&\pi_+(U)e_n=e_{n+1}, \quad n\in\Z,\\
&\pi_-(U)e_n=e_{n+1},\quad n\in\Z\setminus \{-1,0\},\quad
\pi_-(U)e_{-1}=e_1, \quad
\pi_-(U)e_{0}=0.
\end{align*}
Note that the representation $\pi_-$ is non-unital. More precisely, 
$\pi_-(1)$ is the projection onto $\lin\{e_n\,:\, n\in\Z\setminus\{0\}\}$.
Composing $\pi_{\pm}$ with the symbol map $\sigma$
yields the second generator of $\KN(\Cpodl)$, say $[(\epsilon_+,\epsilon_-)]$,
where $\epsilon_\pm = \pi_\pm\circ\hs$.

\subsection{Hopf fibrations of quantum SU(2) over the generic Podle\'s spheres}
                                                                  \label{Hopf}
Let again $q\in(0,1)$.
The Hopf *-algebra $\SUq$ of polynomial functions on the quantum group
$\mathrm{SU}_q(2)$
is generated by $\ha$, $\hb$, $\hc$, $\hd$ with relations
\begin{align*}
& \ha \hb =q \hb \ha,\quad \ha \hc =q\hc \ha,\quad \hb \hd = q \hd \hb,\quad
\hc \hd = q \hd \hc, \quad \hb \hc =\hc \hb,\\
& \ha \hd - q \hb \hc = 1, \quad  \hd \ha - q^{-1} \hb \hc = 1,
\end{align*}
and involution $\ha^*=\hd$, $\hb^*=-q\hc$. The Hopf algebra structure is
given by
\begin{align*}                                       
& \Delta(\ha) = \ha \otimes \ha + \hb \otimes \hc, \quad
\Delta(\hb) = \ha \otimes \hb + \hb \otimes \hd,\\
& \Delta(\hc) = \hc \otimes \ha + \hd \otimes \hc, \quad
\Delta(\hd) = \hc \otimes \hb + \hd \otimes \hd,     \\
& \varepsilon(\ha)=\varepsilon(\hd)=1, \quad
\varepsilon(\hb)=\varepsilon(\hc)=0,\\
& \kappa (\ha)=\hd,\quad \kappa (\hb)=-q^{-1}\hb,
\quad \kappa (\hc)=-q \hc,\quad
\kappa (\hd)=\ha.
\end{align*}

The polynomial *-algebra $\podl$ can be embedded into $\SUq$ by setting
\[                                                   \label{emb}
 \eta_s:=(\hd+q^{-1}s\hb)(\hb-s\hd),\quad \zeta_s:=1-(\ha-qs\hc)(\hd+s\hb).
\]
By restricting the comultiplication $\Delta$ of $\SUq$ to $\podl$,
the latter becomes a right $\SUq$-comodule *-algebra.

Let $J_s$ denote the coalgebra $\SUq/\podl^+\SUq$ with quotient map
$\pr_s:\SUq\rightarrow J_s$,
where $\podl^+=\{x\in\podl\,:\,\varepsilon (x)=0\}$.
It has been shown in \cite{MS} that
$J_s$ is spanned by group-like elements $g_N$,
$N\in\Z$.
Set
$$
M_N:=\{x\in \SUq\,:\, \pr_s(x_{(1)}) \otimes x_{(2)}=g_N\otimes x \}.
$$
Then the algebra of co-invariants $M_0$ is isomorphic
to $\podl$ (via the embedding \eqref{emb}),
$M_N$ is a finitely generated projective $\podl$-module,
and $\SUq$ is the direct sum of all $M_N$ \cite{MS}.
Analogous to the classical Hopf fibration over the 2-sphere, the projective
modules $M_N$ are considered as line bundles over the quantum 2-sphere
$\mathrm{S}^2_{qs}$ with winding number $N$.

Descriptions of
idempotents representing $M_N$ for all $N$ and all parameters
$s$ have been given in \cite{HMS,SW}. For the convenience of the reader,
we review briefly the construction of idempotents from \cite{SW}.

For $N\in\Z$, set
\begin{align}                                              \label{uN}
\begin{split}
 & u_N=q^{-N}(\hb-qs\hd)(\hb-q^2s\hd)\cdots(\hb-q^N s\hd),\quad N>0, \\
& u_{N}=(\hd+q^{-1}s\hb)(\hd+q^{-2}s\hb)\cdots(\hd+q^{N} s\hb),\quad N<0,
\end{split}
\end{align}
and $u_0=1$.
From \cite{BM},
we conclude that $\pr_s(u_{N})=g_{N}$ and thus
$u_N\in M_N$.
(In order to apply the results of \cite{BM},
one has to switch from right coalgebras to left coalgebras by using the
$\ast$-al\-ge\-bra automorphism and coalgebra
anti-homomorphism $\theta:\SUq\rightarrow\SUq$ determined by
$\theta (\ha) = \ha$, $\theta (\hd)=\hd$, $\theta(\hb) =-q\hg$,
$\theta(\hg)=-q^{-1} \hb$.)

For $l\in\frac{1}{2}\N_0$ and $j,k=-l,-l+1,\ldots,l$, let $t^l_{jk}$
denote the matrix elements of finite
dimensional unitary corepresentations of $\SUq$.
By the explicit
description of the  matrix elements (see, e.g., \cite[Section 4.2.4]{KS}), one sees
that $u_{\pm 2l}\in \lin\{t^{l}_{j,\pm l}\,:\, j=-l,-l+1,\ldots,l\}$.
In particular, $u_{\pm 2l}$ is a highest weight vector
of an irreducible spin-$l$ corepresentation. That is,
there are elements $v^l_{k,\pm l}$, $k=-l,-l+1,\ldots,l$, such that
$\Delta(v^l_{k,\pm l})=\sum_j v^l_{j,\pm l}\otimes t^{l}_{j,k}$
and $v^l_{l,\pm l}=c_{\pm l}\, u_{\pm 2l}$ with positive constants
$c_{\pm l}$ determined as follows: Since  $t^l_{jk}$ are
matrix elements of unitary corepresentations, we have
$\Delta(\sum_k v^{l*}_{k,\pm l}v^l_{k,\pm l})=(\sum_k v^{l*}_{k,\pm l}v^l_{k,\pm l})\otimes 1$, so the
sum belongs to the trivial corepresentation.
Thus we find $c_{\pm l}\in\R_+$ such that
$\sum_k v^{l*}_{k,\pm l}v^l_{k,\pm l}=1$. 

Next, define
\[                                                    \label{VN}
V_{\pm 2l}:=(v^l_{-l,\pm l}, \ldots, v^l_{l,\pm l})^\trans
\quad \mbox{and} \quad P_{\pm 2l}= V_{\pm 2l}V_{\pm 2l}^*,
\]
where the superscript $\trans$ denote the transpose. As
$$
V_{\pm 2l}^*V_{\pm 2l}=\sum_k v^{l*}_{k,\pm l}v^l_{k,\pm l}=1,
$$
$P_{\pm 2l}$
is a self-adjoint idempotent. It has been shown in
\cite{SW} that the entries of  $P_{\pm 2l}$
belong to $\podl$ and that
$M_{ N}\cong \podl^{|N|+1}P_{N}$,
where $N=\pm 2l\in\Z$.

Let  $\CSU$ denote the
\Cs-completion of $\SUq$ with respect to the universal \Cs-norm.
Clearly, each *-re\-pre\-sen\-ta\-tion of $\SUq$ restricts to one
of $\podl$. On the other hand, the faithful representation
$\rho_-\oplus\rho_+$ of $\podl$ on $\lN\oplus\lN$ can be realized as a
subrepresentation of a *-re\-pre\-sen\-ta\-tion of $\SUq$, for instance,
of the GNS-re\-pre\-sen\-ta\-tion associated with the Haar state
on $\SUq$ (see, e.g., \cite[Section 6]{SW}).
From the definition of the universal \Cs-norm
and the isomorphism
$\Cpodl\ni\xi\mapsto
(\rho_-\ppt{\oplus}\ppt\rho_+)(\xi)\in(\rho_-\ppt{\oplus}\ppt\rho_+)(\Cpodl)$,
it follows that the embedding is an isometry. Thus $\Cpodl$ can be
viewed as a subalgebra of $\CSU$.

\subsection{Two auxiliary lemmas}                     \label{lem}

The first lemma is used to handle  expressions
like $f(\zeta_s)\eta_s =  \eta_s f(q^2\zeta_s)$
and $f(\zeta_s) S =  S f(q^2\zeta_s)$.

\ble                                               \label{ABC}
Given a unital \Cs-algebra $\A$, let
$A, B,C\in\A$ such that $A^*=A$, $B^*=B$ and $AC=CB$. Then
$f(A)C=Cf(B)$ for any continuous function
$f:\spec(A)\cup\spec(B)\rightarrow \C$.
\ele
\begin{proof}
From $AC=CB$, it follows that
$p(A)C=Cp(B)$ for any polynomial
$p$ in one indeterminate. Approximating $f$ uniformly by polynomials
on the compact
set $\spec(A)\cup\spec(B)$  gives the result since the
multiplication on $\A$ is continuous.
\end{proof}

The second lemma summarizes some crucial relations
between the generators of $\podl$ and highest/lowest
weight vectors of $M_N$.
Before stating the lemma, recall the definition of
the highest weight vectors $u_N\in M_N$
from the previous subsection.
Now, for $N\in\Z$, set
\begin{align}                                                   \label{wN}
\begin{split}
& w_N:=(\ha-qs\hc)(\ha-q^2s\hc)\cdots(\ha-q^N s\hc), \quad N>0,\\
& w_{N}:=(-q)^{-N}(\hc+q^{-1}s\ha)(\hc+q^{-2}s\ha)\cdots(\hc+q^{N} s\ha),\quad N<0,
\end{split}
\end{align}
and $w_0=1$.
Applying the coaction $\Delta$ on $u_N$ and
taking only the elements with $\ha$ and $\hc$ in the left tensor factor,
one sees that  $w_N$ is a lowest
weight vector belonging to the same corepresentation as $u_{N}$.
In particular,  $w_{\pm 2l}= d_{\pm 2l}\, v^l_{-l,\pm l}$
for some non-zero real constant $d_{\pm 2l}$,
where $l\in\frac{1}{2}\N_0$.

\begin{lem}                                              \label{comps}
For $n\in\N_0$,
\begin{align}
&u_n\,\zeta_{q^n\mpt s}=\zeta_{s}\,u_n,\quad
w_n\,\zeta_{q^n\mpt s}=q^{2n} \zeta_{s}\,w_n,              \label{uzwz}\\
&u_{-n}\,\zeta_{q^{-n}\mpt s}=q^{-2n}\zeta_{s}\,u_{-n},\quad
w_{-n}\,\zeta_{q^{-n}\mpt s}=\zeta_{s}\,w_{-n},\\
&u_n\,\eta_{q^n\mpt s}=q^n \eta_{s}\,u_n,\quad
u_{-n}\,\eta_{q^{-n}\mpt s}=q^{-n}\eta_{s}\,u_{-n};
\end{align}
\begin{align}                                              \label{uNuN*}
&u_nu_n^*=q^{-2n} (q^2s^2+\zeta_s) (q^4s^2+\zeta_s)   \cdots
(q^{2n}s^2+\zeta_s),\\
& w_nw_n^*=(1-\zeta_s)(1-q^{2}\zeta_s)\cdots
(1-q^{2(n-1)}\zeta_s),\\
&u_{-n} u_{-n}^*=(1-q^{-2}\zeta_s)(1-q^{-4}\zeta_s)\cdots
(1-q^{-2n}\zeta_s),\\
& w_{-n}w_{-n}^*= (s^2+\zeta_s) (q^{-2}s^2+\zeta_s)   \cdots
(q^{-2(n-1)}s^2+\zeta_s);
\end{align}
\[                                                       \label{uNwN*}
u_nw_n^*=q^{n(n-1)/2}\eta_{s}^n,\quad \ \;
u_{-n}w_{-n}^*=q^{-n(n-1)/2}\eta_{s}^n;
\]
\begin{align}                                              \label{uN*uN}
&u_n^*u_n=q^{-2n}
(q^2s^2+\zeta_{q^n\mpt s}) (q^4s^2+\zeta_{q^n\mpt s})   \cdots
(q^{2n}s^2+\zeta_{q^n\mpt s}),\\
& w_n^*w_n=(1-q^{-2}\zeta_{q^n\mpt s})(1-q^{-4}\zeta_{q^n\mpt s})\cdots
(1-q^{-2n}\zeta_{q^n\mpt s}),                                 \label{wN*wN}\\
&u_{-n}^* u_{-n}=(1-\zeta_{q^{-n}\mpt s})(1-q^{2}\zeta_{q^{-n}\mpt s})\cdots
(1-q^{2(n-1)}\zeta_{q^{-n}\mpt s}),\\
& w_{-n}^*w_{-n}=(s^2+\zeta_{q^{-n}\mpt s})
(q^{-2}s^2+\zeta_{q^{-n}\mpt s})
\cdots  (q^{-2(n-1)}s^2+\zeta_{q^{-n}\mpt s}).
\end{align}
\end{lem}
\begin{proof} The lemma is proved by direct computations using
the embedding \eqref{emb}, the commutation relations in $\SUq$,
and induction on $n$.
\end{proof}

\section{Fibre product approach to line bundles over the generic
Podle\'s spheres}
                                                     \label{sec-glue}

The isomorphism \eqref{CS} admits a geometric interpretation as
describing the generic Podle\'s spheres by gluing two quantum discs
along their boundaries. To see this, note that
$\spec(\rho_-(\zeta_s))\subset [-s^{2}q^{2},0]$ and
$\spec(\rho_+(\zeta_s))\subset [0,q^{2}]$. As a consequence,
$(1-\rho_-(\zeta_s))$ and $(s^2\,{+}\,\rho_+(\zeta_s))$
are strictly positive operators with bounded inverses.
Let $z_-$ and $z_+$ be defined by ``stereographic projection'':
\[
z_-:=s^{-1} \rho_-(\eta_s) (1-\rho_-(\zeta_s))^{-1/2}, \quad
z_+:=\rho_+(\eta_s)(s^2+\rho_+(\zeta_s))^{-1/2}.
\]
Then one readily checks that $z_-$ and $z_+$ satisfy the
defining relation \eqref{Dq} of the quantum disc and therefore generate
$\T$, the \Cs-algebra of $\Dq$.
The generators $\rho_\pm(\zeta_s)$ and
$\rho_\pm(\eta_s)$ can be recovered
from $z_\pm$ since
$$
\rho_-(\zeta_s)= (q^{-2}-1)^{-1}s^2(z_-z_-^*-z_-^*z_-),\quad
\rho_+(\zeta_s)= (q^{-2}-1)^{-1}(z_+z_+^*-z_+^*z_+).
$$
Thus we can view the two copies of $\T$ in \eqref{CS} as
algebras of continuous functions on the quantum disc derived from
continuous functions on the northern hemisphere $\zeta_s\geq 0$
and the southern hemisphere $\zeta_s\leq 0$ by
``stereographic projection''.

Furthermore, the isomorphism \eqref{CS} shows that $\Cpodl$ is
obtained as the fibre product of two Toeplitz algebras
by the following pull-back diagram:
\begin{equation}                                   \label{fp-Cpodl}
        \begin{CD}
    {\Cpodl} @ >{\mathrm{pr}_1}>> {\T} @.\\
    @ V{\mathrm{pr}_0} VV @ V{\hs} VV @.\\
    {\T} @ >{\hs} >> {\CS}. @.\\
        \end{CD}
\end{equation}

In the classical case $q=1$, the surjection $\hs:\T\ra\CS$ corresponds
to an embedding of the circle $\mathrm{S}^1$ into the closed disc
$\mathrm{D}^2$, and the relation $\hs(a_0)=\hs(a_1)$ in \eqref{CS}
means that a pair of continuous functions on the northern
and southern hemispheres is identified with a continuous function
on the 2-sphere iff their restrictions to $\mathrm{S}^1$ coincide.
This identification captures precisely the meaning of
gluing two (quantum) discs along their boundaries.

Recall that complex line bundles with winding number $N\in\Z$ over the classical 2-sphere
can be constructed by taking trivial bundles over the northern and southern hemispheres
and gluing them together along the boundary via the map
$U^{N}:\mathrm{S}^1  \ra \mathrm{S}^1$,\,
$U^{N}(\mathrm{e}^{\im\phi})=\mathrm{e}^{\im\phi N}$.

Our aim is to give a non-commutative analogue of this construction.
The fibre product \eqref{fp-Cpodl} tells us that the trivial bundles over the discs
$\mathrm{D}^2$ should be replaced by $\T$, whereas the transition map $U^{N}$
remains the same. This leads to the following pull-back diagram:
\begin{equation*}                                                         
\xymatrix{
& \T \,{\underset{(U^N\sigma ,\sigma)}{\times}}\, \T
\ar[dl]_{\mathrm{pr}_0} \ar[dr]^{\mathrm{pr}_1}& \\
\T \ar[d]_{\sigma} &    &  \T \ar[d]^{\sigma}\\
\CS \ar[rr]_{f\mapsto U^Nf} & & \CS,   }
\end{equation*}
where
\begin{equation}                                                     \label{fpLB}
 \T \times_{(U^N\sigma ,\sigma)} \T
\cong \{(a_0,a_1)\in \T\oplus \T \,:\, U^{N}\hs(a_0)=\hs(a_1)\}.
\end{equation}

Thinking of $\T$ and $\CS$ as left modules over themselves, it follows from the discussion
in the end of Section \ref{sec-fp} that
$\T\times_{(U^N\sigma ,\sigma)} \T$
is a left $\Cpodl$-module.
This can also be seen directly from Equations \eqref{CS} and\eqref{fpLB}.

The next proposition gives a projective module description of
$\T\times_{(U^N\sigma ,\sigma)} \T$.

\bpr                                                         \label{isoMchi}
For $N\in\Z$, let $L_N:=\T\times_{(U^N\sigma ,\sigma)} \T$,
and define
\[                                                       \label{EN}
E_N := (S^{N}S^{*N},1),\ \ N\geq 0,\quad E_N := (1,S^{|N|}S^{*|N|}),\ \ N< 0.
\]
Then the left $\Cpodl$-modules $L_N$ and $\Cpodl E_N$ are
isomorphic.
\epr

\begin{proof}
Clearly, $E_N$ is a  projection in $\Cpodl$ since $\hs(S^{n}S^{*n})=1=\hs(1)$
and $S^{*n}S^{n}=1 $ for all $n\in\N_0$.

Let $N\geq 0$.
From $\hs(S^{*})=U^{*}$, \eqref{CS} and \eqref{fpLB}, it follows that
$(a_0 S^{*N},a_1)$ belongs to $\Cpodl$ for all $(a_0, a_1)\in L_N$.
Consider now the \Cpodl-linear map
$$
\Psi_N\,:\, L_{N}\rightarrow \Cpodl E_N,\quad
\Psi_N((a_0, a_1)) := (a_0 S^{*N}, a_1)E_N.
$$
We claim that $\Psi_N$ is an isomorphism. Note that
$\Psi_N((a_0, a_1))= (a_0 S^{*N}, a_1)$
since, as above, $S^{*N}S^{N}=1$.
Assume that $(a_0, a_1)\in\ker(\Psi_N)$.
Then $0=(a_0 S^{*N},a_1)$. Hence $a_1=0$ and, as $S^{*N}$ is
right invertible, $a_0=0$. Therefore $\ker(\Psi_N)=\{0\}$, so
$\Psi_N$ is injective.

Next, let $(a_0, a_1)\in \Cpodl$.
Then $(a_0S^N, a_1)\in L_N$ by \eqref{CS} and \eqref{fpLB},  and
$\Psi_N((a_0S^N, a_1))=(a_0S^NS^{*N}, a_1)=
(a_0, a_1)E_N$.
Therefore $\Psi_N$ is also surjective.

For $N<0$, one proves analogously that
$$
\Psi_{N}\,:\, L_{N}\rightarrow \Cpodl E_N,\quad
\Psi_{N}((a_0,a_1)) := (a_0,a_1S^{*|N|})E_N,
$$
is an isomorphism.
\end{proof}

\section{Projective module descriptions
of line bundles associated to the generic Hopf fibration}
                                                        \label{pmd}

In Section \ref{Hopf},
we explained that the quantum line bundles $M_N$ from the generic
Hopf fibration of $\mathrm{SU}_q(2)$ are, as left $\podl$-modules,
isomorphic to $\podl^{|N|+1}P_{N}$, where $P_N$ was given in \eqref{VN}.
As the entries of $P_N$ belong to $\podl$, these projections
represent \Kn-classes of $\Cpodl$.
The aim of this section is to prove that $P_N$ and
the 1-dimensional projection $E_N$ from Equation \eqref{EN}
define the same class in \Kn-theory. As an intermediate step,
we first reduce the $(|N|+1)\times (|N|+1)$-projections $P_N$
to Murray-von Neumann equivalent $2\times 2$-projections
by assembling the corners of $P_N$.
These $2\times 2$-projections are quantum analogues of the classical
Bott projections.
After that,
we prove the Murray-von Neumann equivalence of the ``Bott projections''
and $E_N$.
Note that the 1-dimensional projections $E_N$, $N\neq 0$, do not have classical
counterparts since $\mathrm{S}^2$ is connected and therefore
all projections in $\cC(\mathrm{S}^2)$ are trivial.

Let  $l\in \frac{1}{2}\N$.
Recall from the Preliminaries that
$P_{\pm 2l}= V_{\pm 2l}V_{\pm 2l}^*$, where
$V_{\pm 2l}=(v^l_{-l,\pm l}, \ldots, v^l_{l,\pm l})^\trans$,
and that $v^l_{l,\pm l}\sim  u_{\pm 2l}$ and
$v^l_{-l,\pm l}\sim  w_{\pm 2l}$ with
$u_{\pm 2l}$ and $w_{\pm 2l}$ given in
Equations \eqref{uN} and \eqref{wN}, respectively.
By Lemma \ref{comps}, we have, for $n\in\N$,
\begin{align}
\begin{split}                                 \label{uuww}
&\!\!u_{-n}^{*}u_{-n}^{*}+q^{n(n-1)}w_{-n}^{*}w_{-n}^{*}=
\mbox{$\Pi_{k=0}^{n-1}$} (1\!-\!q^{2k}\zeta_{q^{-n}\mpt s})
+\mbox{$\Pi_{k=0}^{n-1}$}(s^2\!+\!q^{2k}\zeta_{q^{-n}\mpt s}),\\
&\!\!q^{-n(n-1)}u_n^{*}u_n^{*}+w_n^{*}w_n^{*}=
\mbox{$\Pi_{k=0}^{n-1}$}(s^2\!+\!q^{2k-2n}\zeta_{q^n\mpt s})
+\mbox{$\Pi_{k=0}^{n-1}$}(1\!-\!q^{2k-2n}\zeta_{q^n\mpt s}).
\end{split}
\end{align}
Here, the constants $q^{\pm n(n-1)}$ were inserted in order to
obtain more symmetric formulas.
Using \eqref{spec}, one easily shows that
$u_{- n}^{*}u_{- n}^{*}+q^{n(n-1)}w_{- n}^{*}w_{- n}^{*}$ and
$q^{-n(n-1)}u_n^{*}u_n^{*}+w_n^{*}w_n^{*}$
are strictly positive operators with bounded inverses
in $\CSU$. Thus we can define the following 2-vectors
with entries in $\CSU$:
\begin{align}
\begin{split}                                               \label{WN}
  W_{- n}&:=(u_{- n}\,,\, -q^{n(n-1)/2}w_{- n})^\trans
(u_{- n}^*u_{- n}+q^{n(n-1)}w_{- n}^*w_{- n})^{-1/2},\\
W_{n}&:=(q^{-n(n-1)/2}u_{n}\,,\, w_{n})^\trans
(q^{-n(n-1)}u_{n}^*u_{n}+w_{n}^*w_{n})^{-1/2}.
\end{split}
\end{align}
Obviously, $W_{\pm n}^* W_{\pm n} = 1$, so that
\[                                              \label{QN}
Q_{\pm n}:=W_{\pm n}\, W_{\pm n}^*
\]
is a self-adjoint idempotent.

\begin{lem}                                                      \label{LEMQN}
Let $n\in\N$. With $t$ an indeterminate, define
\[                                                     \label{fg}
f_{n}(t):=
\mbox{$\Pi_{k=0}^{n-1}$}(1\mpt -\mpt q^{2k}t), \quad
g_{n}(t):=
\mbox{$\Pi_{k=0}^{n-1}$}(s^2\mpt +\mpt q^{2k}t).
\]
Set $\S:=(S,S)\in\Cpodl$, where $S$ denotes the shift operator on $\lN$.
Then
$$
Q_{-n}=
\begin{pmatrix}
\frac{f_{n}(q^{-2n}\zeta_s)}{f_{n}(q^{-2n}\zeta_s)+g_{n}(q^{-2n}\zeta_s)}
& -\frac{\sqrt{f_{n}(q^{-2n}\zeta_s)g_{n}(q^{-2n}\zeta_s)}}
{f_{n}(q^{-2n}\zeta_s)+g_{n}(q^{-2n}\zeta_s)} \S^n \\[6pt]
-\frac{\sqrt{f_{n}(\zeta_s)g_{n}(\zeta_s)}}
{f_{n}(\zeta_s)+g_{n}(\zeta_s)}\S^{*n}
& \frac{g_{n}(\zeta_s)}{f_{n}(\zeta_s)+g_{n}(\zeta_s)}
\end{pmatrix},
$$
$$
Q_n=
\begin{pmatrix}
\frac{g_n(q^{-2n}\zeta_s)}{f_n(q^{-2n}\zeta_s)+g_n(q^{-2n}\zeta_s)}
& \frac{\sqrt{f_n(q^{-2n}\zeta_s)g_n(q^{-2n}\zeta_s)}}
{f_n(q^{-2n}\zeta_s)+g_n(q^{-2n}\zeta_s)} \S^n \\[6pt]
\frac{\sqrt{f_n(\zeta_s)g_n(\zeta_s)}}{f_n(\zeta_s)+g_n(\zeta_s)}\S^{*n}
& \frac{f_n(\zeta_s)}{f_n(\zeta_s)+g_n(\zeta_s)}
\end{pmatrix}.
$$
In particular, the entries of $Q_{\pm n}$ belong to $\Cpodl$.
\end{lem}
\begin{proof}
The entries $(Q_{\pm n})_{ij}$, $i,j=1,2$, of the
$2\times 2$-matrices $Q_{\pm n}$ are easily computed
by using Lemmas \ref{ABC} and \ref{comps}.
For instance,
\begin{align*}
(Q_n)_{12}&=q^{-n(n-1)/2}u_{n}\ppt
(q^{n(n-1)}u_{ n}^*u_{ n}+w_{n}^*w_{ n})^{-1}\ppt w_{n}^*    \\
&= q^{-n(n-1)/2}u_{n}  \big(\mbox{$\Pi_{k=0}^{n-1}$}(s^2+q^{2k-2n}\zeta_{q^n\mpt s})
+\mbox{$\Pi_{k=0}^{n-1}$}(1-q^{2k-2n}\zeta_{q^n\mpt s})\big)^{-1} w_{n}^*  \\
&=q^{-n(n-1)/2}\big(\pink (s^2+q^{2k-2n}\zeta_{s})+
\pink (1-q^{2k-2n}\zeta_{s}) \big)^{-1} u_{n}  w_{n}^* \\
&=\big( f_n(q^{-2n}\zeta_s)+g_n(q^{-2n}\zeta_s)\big)^{-1}\ppt\eta_{s}^{n},
\end{align*}
where we used \eqref{uuww} in the second equality,
\eqref{uzwz} in the third, and \eqref{uNwN*} and \eqref{fg} in the fourth.
From the second relation in \eqref{rel} and the identification
$\eta_{s}\mapsto (\rho_-(\eta_{s}),\rho_+(\eta_{s}))\in\Cpodl$,
it follows that the polar decomposition of $\eta_{s}$ is given by
$\eta_{s}= \S\sqrt{(1-\zeta_s)(s^2+\zeta_s)}$. Since
$\S\zeta_s= q^{-2}\zeta_s\S$, Lemma \ref{ABC} implies that
$\eta_{s}^{n}=\sqrt{f_n(q^{-2n}\zeta_s)\,g_n(q^{-2n}\zeta_s)}\,\S^{n}$.
Inserting the last equation into the above expression for $ (Q_n)_{12}$
gives the result.

Clearly, $\S^{n}$ and $\S^{*n}$ belong to $\Cpodl$.
Thus, to prove the last statement of the
lemma, it suffices to show that the
functions in $\zeta_{s}$ are continuous on $\spec(\zeta_s)$.
This can easily be verified by observing that
the denominators are non-zero.
\end{proof}

\begin{rem}
Since $u_{\pm n}  w_{\pm n}^*\sim \eta_{s}^{n}$,
the entries of $Q_{\pm n}$ are actually rational functions
in the generators of $\podl$.
This is in analogy to the classical Bott projections
(see, e.g., \cite[Section 2.6]{GFV}).
\end{rem}

Recall the definitions of $P_N$ and $E_N$ in
Equations \eqref{VN} and \eqref{EN}, respectively.
For convenience of notation, set $Q_0:=(1\,,0)^{\trans}(1\,,0)$.
In the next proposition, we show that, for fixed $N\in\Z$, the
projections  $P_N$, $E_N$ and $Q_N$ are
Murray-von Neumann equivalent.

\bpr                                                   \label{QEchi}
For all $N\in\Z$, the projections $P_{N}$, $E_{N}$ and $Q_{N}$
belong to the same class in $\Kn(\Cpodl)$.
\epr

\begin{proof}
The case $N=0$ is trivial. For $N\neq 0$, set
\[
X_{ N} := V_{ N}\, W_{ N}^*,
\]
where $V_{ N}$ and $W_{ N}$ are defined in Equations \eqref{VN} and \eqref{WN},
respectively.
Since $V_{ N}^*V_{ N} =1$ and $ W_{ N}^*W_{ N}=1$, we have
$X_{ N}X_{ N}^*= P_{N}$ by \eqref{VN} and
$X_{ N}^*X_{ N}= Q_N$ by \eqref{QN}.
To prove the Murray-von Neumann equivalence of $P_{N}$ and
$Q_N$, it remains to verify that the entries of $X_{ N}$
belong to $\Cpodl$.
To this end, recall that $u_{  2l}\sim  v^l_{l,  l}$ and
$w_{  2l}\sim v^l_{-l,  l}$ for $l\in\frac{1}{2}\N$, and that
$v^l_{j,  l}v^{l*}_{k,  l}\in\podl$  since $P_{  2l}$ belongs to 
$\mathrm{Mat}_{2l+1, 2l+1}(\podl)$.
Let $r_{2l}(\zeta_{q^{ 2l}s}):=
(q^{-l(2l-1)}u_{2l}^*u_{2l}+w_{2l}^*w_{2l})^{-1}$.
Then, by Lemma \ref{ABC}, Equation \eqref{uzwz} and the foregoing,
\begin{align*}
&(X_{2l})_{j,1}\,\sim \,v^l_{j,l}\,r_{2l}(\zeta_{q^{2l}s})\,u_{2l}^*\,
\sim\, v^l_{j,l}\,v^{l*}_{l, l}\,r_{2l}(\zeta_{s})\,\in\,\Cpodl,\\
&(X_{2l})_{j,2}\,\sim \,v^l_{j,l}\,r_{2l}(\zeta_{q^{2l}s})\,w_{2l}^*\,
\sim\, v^l_{j,l}\,v^{l*}_{-l, l}\,r_{2l}(q^{4l}\zeta_{s})\,\in\,\Cpodl,
\end{align*}
hence $X_{N}\in\mathrm{Mat}_{N+1, 2}(\Cpodl)$ for all $N>0$. A similar
argumentation yields $X_{N}\in\mathrm{Mat}_{|N|+1, 2}(\Cpodl)$
for $N<0$.

Next, for $n\in\N$, define
$$
Y_n:=\left(\mbox{$\frac{\sqrt{g_n(q^{-2n}\zeta_s)}}{\sqrt{f_n(q^{-2n}\zeta_s)+g_n(q^{-2n}\zeta_s)}}\,,\,
 \frac{\sqrt{f_n(q^{-2n}\zeta_s)}}
{\sqrt{f_n(q^{-2n}\zeta_s)+g_n(q^{-2n}\zeta_s)}}$} \S^n \right),
$$
$$
Y_{-n}:=\left(\mbox{$ \frac{\sqrt{f_{n}(q^{-2n}\zeta_s)}}
{\sqrt{f_{n}(q^{-2n}\zeta_s)+g_{n}(q^{-2n}\zeta_s)}} \,,\,
\frac{\sqrt{g_{n}(q^{-2n}\zeta_s)}}
{\sqrt{f_{n}(q^{-2n}\zeta_s)+g_{n}(q^{-2n}\zeta_s)}}
$} \S^n \right).
$$
By the same arguments as in the proof of Lemma \ref{LEMQN},
it follows that the entries of $Y_{\pm n}$ belong to $\Cpodl$.
Since, by Lemma \ref{ABC},  $\S^{*} f(\zeta_s)=f(q^2\zeta_s)\S^{*}$ for all
continuous functions $f$ on $\spec(\zeta_s)$,
and since $\S^{*n}\S^{n}=1$,
one sees immediately that $Y_{\pm n}^{*}Y_{\pm n}=Q_{\pm n}$.
Thus, the proposition will be proved, if we show that
$Y_{\pm n}Y_{\pm n}^{*}=E_{\pm n}$.

Observe that
$\rho_-(g_{ n}(q^{-2n}\zeta_s))e_k=0$ and
$\rho_+(f_{n}(q^{-2n}\zeta_s))e_k=0$ for all
$k=0,\ldots, n-1$,
and that
\begin{align*}
&\rho_-(f_{n}(q^{-2n}\zeta_s))e_j
=f_{n}(-s^{2}q^{2(j+1-n)})e_j \neq 0,\\
&\rho_+(g_{n}(q^{-2n}\zeta_s))e_j
=g_{n}(q^{2(j+1-n)})e_j\neq 0, \qquad j\in\N_0.
\end{align*}

Using the identification
$h\mapsto (\rho_-(h),\rho_+(h))$
for elements $h\in\Cpodl$,
one concludes that
\begin{align}
\begin{split}                                                  \label{SS}
 \mbox{$\frac{g_n(q^{-2n}\zeta_s)}{f_n(q^{-2n}\zeta_s)+g_n(q^{-2n}\zeta_s)}$}
(1\!-\!S^{n}S^{*n}\,,\,1\!-\!S^{n}S^{*n}) &=(0\,,\,1\!-\!S^{n}S^{*n}),\\
\mbox{$\frac{f_n(q^{-2n}\zeta_s)}{f_n(q^{-2n}\zeta_s)+g_n(q^{-2n}\zeta_s)}$}
(1\!-\!S^{n}S^{*n}\,,\,1\!-\!S^{n}S^{*n}) &=(1\!-\!S^{n}S^{*n}\,,\,0).
\end{split}
\end{align}
Since $\S^{n}\S^{*n}=(S^{n}S^{*n}\,,\,S^{n}S^{*n})$, we get
\begin{align*}
 Y_{ n}Y_{n}^{*}&=\Big(
\mbox{$\frac{g_n(q^{-2n}\zeta_s)}{f_n(q^{-2n}\zeta_s)+g_n(q^{-2n}\zeta_s)}
+\frac{f_n(q^{-2n}\zeta_s)}{f_n(q^{-2n}\zeta_s)+g_n(q^{-2n}\zeta_s)}$}\Big)\S^{n}\S^{*n}\\
&\qquad\qquad\qquad\qquad\qquad\qquad\qquad
+\mbox{$\frac{g_n(q^{-2n}\zeta_s)}{f_n(q^{-2n}\zeta_s)+g_n(q^{-2n}\zeta_s)}$}(1-\S^{n}\S^{*n})\\
&=(S^{n}S^{*n}\,,\,S^{n}S^{*n})+(0\,,\,1\!-\!S^{n}S^{*n})= E_{n},
\end{align*}
and similarly $Y_{-n}Y_{-n}^{*}= E_{-n}$. This concludes the proof.
\end{proof}

From Proposition \ref{isoMchi} and the proof of Proposition \ref{QEchi},
we obtain the following isomorphisms of left $\Cpodl$-modules:
\begin{equation*}
 L_N:=\T\mpt\times_{(U^N\sigma ,\sigma)}\mpt\mpt \T
\cong \Cpodl E_N\cong \Cpodl^{2}Q_N\cong \Cpodl^{|N|+1}P_N,\quad N\in\Z.
\end{equation*}
The common $\Kn$-class will be denoted by $[L_N]$.

Classically, the direct sum of two line bundles over the 2-sphere
with winding number $N$ and $-N$
is isomorphic to a trivial rank 2 bundle.
An analogous result holds in the quantum case:

\bco
The direct sum of the line bundles $L_N$ and $L_{-N}$ is
isomorphic to a free rank 2 module.
\eco
\begin{proof} For $N\neq 0$,
we have obviously $Q_N+Q_{-N}=1\in \mathrm{Mat}_{2, 2}(\Cpodl)$, so
$Q_N$ and  $Q_{-N}= 1-Q_N$ are complementary projections in
$\mathrm{Mat}_{2, 2}(\Cpodl)$. Therefore
$L_N\oplus L_{-N}\cong \Cpodl^{2}Q_N\oplus \Cpodl^{2}Q_{-N}\cong \Cpodl^{2}$.
The case $N=0$ is trivial.
\end{proof}

\section{Index computation for quantum line bundles}


Let $A$ be a C*-algebra, $p\in\Mat_{n,n}(A)$ a projection,
and $\varrho_+$ and $\varrho_-$ *-re\-pre\-sen\-ta\-tions
of $A$ on a Hilbert space $\H$ such that
$[(\varrho_+,\varrho_-)]\in K^{0}(A)$.
If the following traces exist, then the formula
\begin{equation}                                                    \label{CCpair}
\langle [(\varrho_+,\varrho_-)],[p]\rangle=\tr_\H (\tr_{\Mat_{n,n}} (\varrho_+-\varrho_-)(p))
\end{equation}
computes the index of the Fredholm operator
$\varrho_+(p)\varrho_-(p):  \varrho_-(p)\H^{n} \rightarrow \varrho_+(p)\H^{n}$ and
therefore yields a pairing between $K^{0}(A)$ and $K_{0}(A)$
\cite[Section IV.1]{C}.

In general, the difficulty of computing the traces increases with
the growing size of the matrices.
However, the fibre product approach to quantum line bundles
provided us with 1-dimensional projections as representatives
of \Kn-classes. Their simple form makes them very suitable for the
calculation of the index pairing.

To compute the Fredholm indices (Chern numbers) for
quantum line bundles, we shall use the generators $[(\epsilon_+,\epsilon_-)]$
and $[(\rho_+,\rho_-)]$ of $\KN(\Cpodl)$ defined in
the end of Section \ref{Ps}.

\bpr
The pairing between the generators
of $\KN(\Cpodl)$
and the $\Kn$-class of the quantum line bundle $L_N$, $N\in\Z$,
yields
$$
\<\, [(\epsilon_+,\epsilon_-)] \,, [L_N]\,\>=1, \qquad
\<\, [(\rho_+,\rho_-)]  \,, [L_N]\,\>= N.
$$
\epr
\begin{proof}

As already mentioned,  we choose the 1-dimensional projection $E_N$
as a representative of the $\Kn$-class $[L_N]$.
Recall from Section \ref{Ps} that $\pi_+(1)=1$ and
$\pi_-(1)$ is the projection onto $\lin\{e_k\,:\, k\in\Z\setminus\{0\}\}$.
Thus, $(\pi_+-\pi_-)(1)$ is the projection onto $\C e_0$.
Moreover, $\hs((S^{n}S^{*n},1))=\hs((1,S^{n}S^{*n}))=1$ for all $n\in\N_0$,
therefore $\hs(E_N)=1$. Since $\epsilon_\pm=\pi_\pm\circ\hs$, we get
$$
\<\, [(\epsilon_+,\epsilon_-)] \,, [L_N]\,\>
=\tr_{{\ell}^2(\Z)} (\pi_+-\pi_-)(\hs(E_N))
=\tr_{{\ell}^2(\Z)} (\pi_+-\pi_-)(1) =1.
$$

Let $N\geq 0$. By \eqref{rhopm},
$(\rho_+-\rho_-)(E_N)=(\rho_+-\rho_-)(S^{N}S^{*N},1)=1-S^{N}S^{*N}$
is the projection onto $\lin\{e_0,\ldots, e_{N-1}\}$. Thus
$$
\<\, [(\rho_+,\rho_-)]  \,, [L_N]\,\>=\tr_{{\ell}^2(\N_0)} (\rho_+-\rho_-)(E_N)
=\tr_{{\ell}^2(\N_0)}(1-S^{N}S^{*N})=N.
$$
If $N<0$, we get
$(\rho_+-\rho_-)(E_N)=(\rho_+-\rho_-)(1,S^{|N|}S^{*|N|})=S^{|N|}S^{*|N|}-1$,
again by \eqref{rhopm}. Hence
$$
\<\, [(\rho_+,\rho_-)]  \,, [L_N]\,\>=\tr_{{\ell}^2(\N_0)} (\rho_+-\rho_-)(E_N)
=-\tr_{{\ell}^2(\N_0)}(1-S^{|N|}S^{*|N|})=-|N|.
$$
Therefore $\<\, [(\rho_+,\rho_-)]  \,, [L_N]\,\>=N$ for all $N\in\Z$.
\end{proof}

We remark that
the additive map
$\<\,[(\epsilon_+,\epsilon_-)]\,, [\,\cdot\,]\,\>:\Kn(\Cpodl)\rightarrow \Z$
can be used to detect the rank of the bundle (equal to 1 for line bundles).
A similar statement was made in
\cite[Remark 3.4]{HM} for the standard Podle\'s sphere.

The Chern number $\<\, [(\rho_+,\rho_-)]  \,, [L_N]\,\>= N$
coincides with the power of $U$ in \eqref{fpLB}
and thus computes the ``winding number'', i.e.,
the number of rotations along the equator of the two glued
quantum discs.
That the pairing between $[(\rho_+,\rho_-)] $
and the \Kn-class of $P_N$ gives the winding number $N$
has been shown before in \cite{HMS}, 
relying heavily on the index theorem.
The advantage of our fibre product approach is that we are
able to compute the index pairing directly
by finding projections of small matrix size.

\section*{Acknowledgments}
%
The author thanks Piotr Hajac for very useful discussions on the subject.
This work was
supported by the DFG fellowship WA 1698/2-1 and the European Commission grant
MTK-CT-2004-509794.



\begin{thebibliography}{32}
\bibitem{Ba}
Bass,~H.: 
Algebraic K-theory, Benjamin, New York, 1968.

\bibitem{BHMS}
Baum,~P.~F., P.~M.~Hajac, R.~Matthes, W.~Szyma\'nski: 
\emph{The $K$-theory of Heegaard-type
quantum 3-spheres.} $K$-Theory {\bf 35} (2005),  159--186.

\bibitem{Brz}
Brzezi\'nski,~T.:
\emph{Quantum Homogeneous Spaces as Quantum Quotient Spaces.} 
J. Math. Phys. {\bf 37} (1996), 2388--2399.

\bibitem{BM}
Brzezi\'nski,~T., S.~Majid:
\emph{Quantum geometry of algebra factorisations and coalgebra bundles.}
Commun. Math. Phys. {\bf 213} (2000), 491--521.

\bibitem{bk96} Budzy\'nski,~R.~J., W.~Kondracki: 
\emph{Quantum principal fibre bundles: Topological aspects.}
Rep. Math. Phys. {\bf 37} (1996), 365--385.

\bibitem{cm00} Calow,~D., R.~Matthes: 
\emph{Covering and gluing of algebras and differential algebras.} 
J. Geom. Phys. {\bf 32} (2000), 364--396. 

\bibitem{C}
Connes,~A.: Noncommutative geometry.
Academic Press, San Diego, 1994.

\bibitem{DHHM}
D\c{a}browski,~L., T.~Hadfield, P.~M.~Hajac, R.~Matthes:
\emph{K-theoretic construction of noncommutative instantons of all charges.} 
Preprint, arXiv: math/0702001v1.

\bibitem{GFV}
Gracia-Bond\'{\i}a,~J.~M., H.~Figueroa, J.~C.~V\'arilly:
Elements of Noncommutative Geometry.
Birkh\"auser, Boston, 2001.

\bibitem{H}
Hajac,~P.~M.:
\emph{Bundles over quantum sphere and noncommutative index theorem.}
$K$-Theory {\bf 21} (1996), 141--150.

\bibitem{HM}
Hajac,~P.~M., S.~Majid:
\emph{Projective module description of the $q$-mono\-pole.}
Commun. Math. Phys. {\bf 206} (1999), 247--264.

\bibitem{HMS}
Hajac,~P.~M., R.~Matthes, W.~Szyma\'nski:
\emph{Chern numbers for two families of noncommutative Hopf fibrations.}
C. R. Math. Acad. Sci. Paris  {\bf 336}  (2003), 925--930.

\bibitem{HMS06}
Hajac,~P.~M., R.~Matthes, W.~Szyma\'nski:
\emph{Noncommutative index theory for mirror quantum spheres.}
C. R. Math. Acad. Sci. Paris  {\bf 343}  (2006), 731--736.

\bibitem{HW}
Hajac,~P.~M., E.~Wagner:
\emph{The pullbacks of principal coactions.}
In preparation.

\bibitem{kl93}  Klimek~S., A.~Lesniewski:
\emph{A two-parameter quantum deformation of the unit disc.} 
J.\ Funct.\ Anal.\ {\bf 115} (1993), 1--23.

\bibitem{MNW1}
Masuda,~T., Y.~Nakagami, J.~Watanabe:
\emph{Noncommutative
differential geometry on the quantum $\mathrm{SU}(2)$.
I: An Algebraic Viewpoint.}
$K$-Theory {\bf 4} (1990), 157--180.

\bibitem{MNW}
Masuda,~T., Y.~Nakagami, J.~Watanabe:
\emph{Noncommutative
differential geometry on the quantum two sphere of Podle\'s.
I: An Algebraic Viewpoint.}
$K$-Theory {\bf 5} (1991), 151--175.

\bibitem{MS}
M\"uller,~E.~F., H.-J.~Schneider:
\emph{Quantum homogeneous spaces with faithfully flat module structures.}
Israel J. Math. {\bf 111} (1999), 157--190.

\bibitem{KS}
Klimyk,~K.~A., K.~Schm\"udgen: 
Quantum Groups and Their Representations. Springer, Berlin, 1997.

\bibitem{Podles}
Podle\'s,~P.:
\emph{Quantum spheres.}
Lett. Math. Phys. {\bf 14}  (1987), 193--202.


\bibitem{SW}
Schm\"udgen,~K., E.~Wagner:
\emph{Representations of cross product algebras of Podle\'s quantum spheres.}
J. Lie Theory {\bf 17} (2007), 751--790.

\bibitem{s-a91} Sheu,~A.~J.-L.:
\emph{Quantization of the Poisson $SU(2)$ and its Poisson homogeneous space
-- the 2-sphere.}
Commun.\ Math.\ Phys.\ {\bf 135} (1991), 217--232.
\end{thebibliography}
\end{document}